\documentclass{amsproc}
\usepackage[mathcal]{euscript}

\newtheorem{theorem}{Theorem}[section]
\newtheorem{lemma}[theorem]{Lemma}
\newtheorem{corollary}[theorem]{Corollary}

\theoremstyle{definition}
\newtheorem{definition}[theorem]{Definition}

\theoremstyle{remark}
\newtheorem{remark}[theorem]{Remark}

\numberwithin{equation}{section}

\newcommand{\Q}{\begin{pmatrix} e & a & x \\
                     a^* & e & b\\
                     x^* & b^* & e\\
      \end{pmatrix}}
\newcommand{\apro}{\begin{pmatrix} e & a  \\
                     a^* & e \\
      \end{pmatrix}}
\newcommand{\bpro}{\begin{pmatrix} e & b  \\
                     b^* & e \\
      \end{pmatrix}}
      \newcommand{\abpro}{\begin{pmatrix} e & ab  \\
                     (ab)^* & e \\
      \end{pmatrix}}
\newcommand{\Qp}{\begin{pmatrix} 1 & p(a) & p(x) \\
                     \overline{p(a)} & 1 & p(b)\\
                     \overline{p(x)} & \overline{p(b)}  & 1\\
      \end{pmatrix}}

\newcommand{\mna}{M_n(\EuScript{A})}

\newcommand{\spa}{\hspace{0.5em}}
\newcommand{\mnc}{M_n(\mathbb{C})}
\newcommand{\lj}{\lambda_j}
\newcommand{\lk}{\lambda_k}
\newcommand{\lam}{\lambda}

\newcommand{\xab}{e^{i\phi_k}(x_{jk}-\lj
b_{jk})+e^{-i\phi_k}(x_{kj}-\lj b_{kj})}
\newcommand{\QED}{\hfill $\Box$}
\newcommand{\pf}{\noindent{\textbf{Proof: }}}

\newcommand{\real}{\operatorname{Re}}

\newcommand{\ah}{\hat{a}}
\newcommand{\bh}{\hat{b}}
\newcommand{\xh}{\hat{x}}
\renewcommand{\AA}{\EuScript{A}}
\newcommand{\CC}{\mathbb{C}}
\newcommand{\cs}{C^*}

\newcommand{\ZZ}{\mathbb{Z}}
\newcommand{\RR}{\mathbb{R}}
\newcommand{\NN}{\mathbb{N}}
\newcommand{\HH}{\EuScript{H}}

\begin{document}

\title{An Explicit Duality for Finite Groups}

%    Information for first Author
\author{Robert A. Cohen}
\address{Department of Mathematics, University of Colorado,
Campus Box 395, Boulder, Colorado, 80309-0395}

\email{robert.cohen@euclid.colorado.edu}

\author{Martin E. Walter}
\address{Department of Mathematics, University of Colorado,
Campus Box 395, Boulder, Colorado, 80309-0395}
\email{martin.walter@euclid.colorado.edu}

%\dedicatory{This paper is dedicated to .....}


\begin{abstract}
Using a ``3 by 3 matrix trick'' we previously showed that
multiplication in a $C$*-algebra ${\AA}$, an algebraic structure,
is determined by the geometry of the $C$*-algebra of the 3 by 3
matrices with entries from ${\AA}$, $M_3 ({\AA})$.  As an
application of this algebra-geometry duality we now construct an
order theoretic based duality theory for all groups which are
either locally compact abelian or finite.  This construction
generalizes the van Kampen--Pontriagin duality for locally compact
abelian groups. 
\end{abstract}

\maketitle

\specialsection*{Introduction}
 For any locally compact group $G$, the convex, partially ordered
 semigroup $P(G)$ of
 continuous positive definite functions on $G$, is a complete
 invariant of the group.  In the discussion below, we show how
 $P(G)$ can be
 used to recover the algebraic structure of $G$ when $G$ is finite
 or abelian.  In the process we outline a duality theory for
 locally compact groups

Throughout the paper, $\cs (G)$ and $L^1(G)$ will be defined as in
\cite{D}, 13.9.1.  For a $\cs$-algebra $\AA$, we will use $\AA^+$
to denote the positive part of $\AA$, also as in \cite{D}.  The
identity operator of $\AA$ will be denoted by $I$ and $B(\HH)$
stands for the bounded linear operators on the Hilbert space
$\HH$.  For $a,b\in \HH$, the convex hull of $a$ and $b$ is
written as co($a,b$) and the orthogonal complement of the vector
$\xi\in\HH$ is written as $\xi^{\bot}$. The notation $M_n(\AA)$
denotes the $n$ by $n$ matrices with entries from $\AA$; e.g.,
$\AA=\cs (G)$ and $\AA=\CC$ are important examples used in this
paper. The notation diag$(\lam_1,...,\lam_n)$ will be used to
denote a diagonal $n$ by $n$ matrix with $\lam_k$ in the $k^{th}$
diagonal entry. The symbol $\ZZ_n$ will be used to denote the
integers modulo $n$.

\section{Binary Product and Order Structure Duality}
For any locally compact group $G$, the product structure of $G$ is
determined by the order structure of $M_3(C^*(G))$.  The
correspondence can be seen as a special case of the following
theorem:
\begin{theorem}\label{martysthm} Let $\AA$ be a $C^*$-algebra with unit $e$.
Suppose $a,b$ are unitary elements of $\AA$, i.e.,
$a^*a=aa^*=b^*b=bb^*=e$. Then if $x\in\AA$ we have:
\[
   \begin{pmatrix} e & a & x \\
                     a^* & e & b\\
                     x^* & b^* & e\\
      \end{pmatrix}\in M_3(\AA)^+ \]
      if and only if $x=ab$.
\end{theorem}
\proof See \cite{W5}. \QED

If the hypothesis of Theorem \ref{martysthm}
is true and $Q:=\begin{pmatrix} e & a & x \\
                     a^* & e & b\\
                     x^* & b^* & e\\
      \end{pmatrix}\in M_3(\AA)^+$, then $\frac{1}{3}Q$ is a
      projection in $M_3(C^*(G))$.  Unitaries and projections play
      such significant roles in operator algebra theory that an
      interesting observation must be made.  We traditionally
      identify a group element $a$ with a unitary operator in a
      $\cs$-algebra.  In the case of finite groups, we may
      identify $a$ with a unitary matrix with entries from the
      complex numbers (this fact is
      useful for many results in this paper).  However
      by canonically identifying the unitary operator $a$ with the
      block matrix $P=\frac{1}{2}\apro$ we see that we can also identify
      the group element $a$ with a projection operator, $P$.
      Theorem \ref{martysthm}
      shows further that we can construct the product $ab$ via
      concatenation of $\apro$ with $\bpro$
      to produce the positive matrix
      $\begin{pmatrix} e & a & ab \\
                     a^* & e & b\\
                     (ab)^* & b^* & e\\
      \end{pmatrix}$, where the product $ab$ corresponds to the
      two by two submatrix $\abpro$.
\section{Constructing a Dual of $P_1(G)$}

 We often consider the semigroup
$P_1(G)$ consisting of positive definite functions on $G$
satisfying $p(e)\leq 1$.  Indeed, $P_1(G)$ also suffices as a
complete invariant of the group $G$ \cite{W4}.  Our goal is to
define the ``dual of $P_1 (G)$.'' To this end we first consider
certain morphisms of $P_1 (G)$. Let $P(\Gamma_n )$ be defined to
be $\{ M_n ({\Bbb C})^+, \circ \}$, the $n$ by $n$ complex,
positive definite matrices with Schur--Hadamard product, i.e., if
$(p_{jk} ), (q_{jk} )$ are in $M_n ({\Bbb C} )^+$, then $(p_{jk} )
\circ (q_{jk} ) =  (p_{jk} q_{jk} )$, the componentwise product.
Note that this collection of matrices happens to be identifiable
with the set of complex positive definite functions on the
groupoid $\Gamma_n$, of $n$ by $n$ matrix units.  Hence we
introduce the notation $P(\Gamma_n )$, cf., \cite{W3}. Given $n =
1,2,\dots$ consider maps $\varphi_n : P_1 (G) \rightarrow
P(\Gamma_n )$ which satisfy: \item {(I)} $\varphi_n (\lambda p
+(1- \lambda) q)) = \lambda \varphi_n (p) + (1-\lambda )\varphi_n
(q)$, $p,q \in P_1 (G)$, $0 \leq \lambda \leq 1$, i.e.,
$\varphi_n$ is affine; \item {(II)} $\varphi_n (pq) = \varphi_n
(p) \circ \varphi_n (q)$, $p,q \in P_1 (G)$, i.e., $\varphi_n$ is
multiplicative; \item {(III)} For $p,q \in P_1 (G)$, if $p\geq q$,
then $\varphi_n (p) \geq \varphi_n (q)$, i.e., $\varphi_n$ is
order-preserving.

What form must a matrix $\varphi_n (p)$ have?  From the theory of
the Fourier--Stieltjes algebra, $B(G)$, of $G$, cf., \cite{E}
there are elements $x_{jk}$ in the spectrum of this commutative
Banach algebra such that $\varphi_n (p)_{jk} =p(x_{jk} )$ for all
$p \in P_1 (G)$. There are two possibilities for $x_{jk}$.  The
first is that $x_{jk}$ is non-singular, i.e., there is a $g_{jk}
\in G$ such that $p(x_{jk} ) = p(g_{jk} )$ for all $p \in P_1
(G)$.  The second possibility is that $p(x_{jk} ) = 0$ for all $p
\in P_1 (G)$ with compact support.  Such $x_{jk}$ are called {\it
singular\/}.  We thus impose a fourth condition. \item {(IV)} For
each choice of $j$ and $k$, $1 \leq j,k \leq n$, $p \in P_1 (G)
\mapsto \varphi_n (p)_{jk} \in {\Bbb C}$ is nonvanishing on at
least one $p \in P_1 (G)$ with compact support.

\begin{definition}\label{nmorph}A map $\varphi_n : P_1 (G) \rightarrow P( \Gamma_n )$
satisfying I, II, III, IV is called a non-singular $n$-morphism of
$P_1 (G)$, i.e., $\varphi_n$ is affine, multiplicative, order
preserving, and nonsingular.\end{definition}

\begin{remark} For fixed $n = 1,2,\dots ,$ a non-singular
$n$-morphism of $P_1 (G)$ is of the form $p \in P_1 (G) \mapsto
(p(g_{jk} ))_{1\leq j,k \leq n} \in P(\Gamma_n )$ for some $g_{jk}
\in G$, $1 \leq j,k \leq n$.
\end{remark}

Using the enveloping $C$*-algebra and $W$*-algebra of $G$, cf.,
\cite{W1} we can consider $G$ as a collection of unitaries in a
$C$*-algebra and apply Theorem \ref{martysthm} to refine the above
remark, to get

\begin{remark} If $\varphi_n : P_1 (G) \rightarrow
P(\Gamma_n )$ is a non-singular $n$-morphism, then we hope in
general that there exists $g_1 ,g_2 ,\dots ,g_n \in G$ such that
\[\varphi_n (p) = (p(g_j^{-1} g_k))_{1\leq j,k\leq n} \in
P(\Gamma_n ) \text{ for all } p \in P_1 (G).\]  This hope is
realized below in the case when $G$ is abelian or finite.
\end{remark}

We can now define the dual of $P_1 (G)$.

\begin{definition}\label{pdaul} The dual of $P_1 (G)$ consists of those
functions $f : P_1 (G) \mapsto {\Bbb C}$ such that there exists a
non-singular $n$-morphism, $\varphi_n$, satisfying $f(p) =
\varphi_n (p)_{jk}$ for some fixed $j,k$, $1 \leq j, k \leq n$,
and all $p \in P_1 (G)$.  As we will see, the $\varphi_n$ can be
seen as ``multiplication
 table" morphisms.\end{definition}

  It is easy to see that the dual of $P_1(G)$ as a topological
  space is homeomorphic to $G$, but how does one discover its group
structure?  Pick any non-singular $n$-morphism $\varphi_n$ and
suppose $f$ and $h$ are in the dual of $P_1 (G)$ and $f(p) =
\varphi_n (p)_{j_0 k_0}$, $h(p) = \varphi_n (p)_{r_0 s_0}$ for all
$p \in P_1 (G)$.  Thus we may have the following situation if $j_0
< k_0$, $r_0 < s_0$, $j_0 < r_0$, $k_0 < s_0$:
$$\varphi_n (p) = \begin{pmatrix}  p(e) &&&&&&&\\
&\ddots &&&&&&\\
&&p(e) &\cdots &f(p) &\cdots &fh(p) &\\
&&&\ddots &&&&\\
&&&&p(e) &\cdots &h(p) &\\
&&&&&\ddots &\vdots &\\
&&&&&&\ddots &\\
&&&&&&&p(e)\\
&&&&k_0 &&s_0 &\end{pmatrix} \begin{matrix} \phantom{0}\\ \phantom{0}\\
j_0^{\phantom
{0}}\\
\phantom{0_0^0}\\ r_0 \\ \phantom{0}\\ \phantom{0}\\ \phantom{0}\\
\phantom{0}\end{matrix}$$ The product $fh$ should appear as the
element in the $j_0, s_0^{th}$ entry for $\varphi_n$, i.e.,
$(fh)(p) = \varphi_n (p)_{j_0 s_0}$ for all $p \in P_1 (G)$.

This product will be well defined modulo a detail discussed in the
next section and (by checking diagrams) associative.  The identity
is any diagonal entry of any non-singular $n$-morphism, and
$f^{-1}$ is found to satisfy $f^{-1} (p) = \varphi_n (p)_{kj}$
whenever $f(p) = \varphi_n (p)_{jk}$ for particular, $1 \leq j,k
\leq n$, and all $p \in P_1 (G)$.  In the rest of this paper we
outline rigorous proofs of the above construction when $G$ is
abelian or finite.

\section{Recovering Product Structure of $G$ with $P_1(G)$}
 As noted previously, we have that $P_1(G)$ is a complete
invariant of the group $G$.  It follows that, given group elements
$a$ and $b$, $P_1(G)$ can sufficiently determine the product of
$a$ and $b$, but how explicitly can we construct this product?
When trying to answer this question we were led to consider the
elements of $P_1(G)$ in conjunction with Theorem \ref{martysthm}.
Translating an exercise in \cite{K}, we establish a motivating
lemma:
\begin{lemma}\label{kad}
If $[a_{jk}]\in \mna^+$ where $\AA$ is a unital $\cs$-algebra and
$p$ is any state on $\AA$, i.e., $p$ is positive definite and
$p(e)=1$, then $[p(a_{jk})]\in\mnc^+$.\end{lemma}

\proof See \cite{K}, p.884 . \QED

\medskip
 The semigroup $P(G)$ identifies with the positive forms on $\cs
(G)$, c.f., \cite{D}, 13.4.   Without losing generality, we
consider the sub-semigroup $K(G)$ of $P_1(G)$ consisting of
positive definite functions on $G$ satisfying $p(e)= 1$.  We call
these functions \textit{states} on $G$ and call $K(G)$ the
\textit{state space} of $G$.  Note that any state on $G$ can be
uniquely extended to a state on $\cs(G)$.  Thus, by Lemma
\ref{kad},
\[\Q \in M_3(\cs(G))^+
\Rightarrow \Qp\in M_3(\CC)^+ \spa\forall\spa p\in K(G).\]

 The converse of the
 implication above is not true in the non-abelian case.
 This can be seen by noting
that, given $a$ and $b$ in non-abelian $G$ with $ab\neq ba$, for
all $p\in K(G)$ we have
\[\begin{pmatrix} 1 & p(a) & p(ab) \\
                     \overline{p(a)} & 1 & p(b)\\
                      \overline{p(ab)} & \overline{p(b)}  & 1\\
                      \end{pmatrix}\in M_3(\CC)^+\spa \text{ and }
                      \spa\begin{pmatrix} 1 & p(a) & p(ba) \\
                     \overline{p(a)} & 1 & p(b)\\
                      \overline{p(ba)} & \overline{p(b)}  & 1\\
                      \end{pmatrix}\in M_3(\CC)^+
                      .\]
Why?  Because any positive definite function on $G$ is positive
definite on $G^{op}$ (the ``opposed" group where the product is
reversed).  The matrices above are positive for all $p\in K(G)$ by
Theorem \ref{martysthm} and Lemma \ref{kad} applied to $\cs (G)$
and $\cs (G^{op})$, respectively.  However,
\[\begin{pmatrix} e & a & ba \\
                     a^* & e & b\\
                     (ba)^* & b^* & e\\
      \end{pmatrix}\notin M_3(\cs (G))^+.\]

While slightly disappointing, the observation just described does
not dash all hopes of constructing a ``dual" version of Theorem
\ref{martysthm}.  In fact, when $G$ is abelian we immediately
arrive at a pleasing and insightful result.  In the following
theorem we also define an important class \{$Q_p$\} of complex
matrices indexed by $p\in K(G)$.

\begin{theorem} \label{abelian}
Let G be a locally compact abelian group with $a,b, x\in G$. Then
\[ Q_p := \Qp \in M_3(\CC)^+ \spa \forall \spa p\in K(G)\]
if and only if $x=ab$.\end{theorem}

\proof Checking that the determinant of $Q_p$ is positive leads to
(after a little algebraic manipulation) a fundamental inequality:
\begin{equation}\label{fundin}
|p(x)-p(a)p(b)|^2 \leq (1-|p(a)|^2)(1-|p(b)|^2) \end{equation}

  If $p$ is any character on $G$, then $|p(a)|=1$.
Thus by \eqref{fundin}, since characters are multiplicative, we
have $p(x)=p(a)p(b)=p(ab)$.  Also, since characters separate the
points of $G$, we have $x=ab$.

   The other direction is
trivial using Theorem \ref{martysthm} and Lemma \ref{kad}.  \QED

\medskip
  Returning to the nonabelian case, it follows from \eqref{fundin}
  that if $|p(a)|=1$ and $Q_p \in M_3(\CC)^+$ then $p(x)=p(a)p(b)$.
  Recalling that $Q_p\in
M_3(\CC)^+$ when $x=ab$ or $x=ba$; we also have by \eqref{fundin}
that if $|p(a)|=1$ for some $p$, then $p(a)p(b)=p(ab)=p(ba)$.
Finally, assuming $Q_p\in M_3(\CC)^+$ for some $x\in G$, we get
that $p(x)=p(a)p(b)=p(ab)=p(ba)$ whenever $|p(a)|=1$; or
equivalently, whenever $|p(b)|=1$. This suggests that $x$ might be
$ab$ or $ba$.
 \begin{definition}  For a unitary element $a$ of a unital $\cs$-algebra
 $\AA$, a state $p$ on $\AA$ such that $|p(a)|=1$ is called an
 \textit{eigenstate} of
$a$. \end{definition}

These eigenstates are fundamental to our study.  Note that the
determinant of $Q_p$ is zero at any eigenstate of $a$.  Also, as
just noted, $x$ agrees with $ab$ and $ba$ at any eigenstate of $a$
or $b$.  In other words, if $a$, $b$, and $x$ are unitary complex
matrices written with respect to a basis of eigenvectors of $a$,
then we immediately see that $x$ agrees with $ab$ and $ba$ on the
diagonal.

We will now state our main result for a ``dual" version of Theorem
\ref{martysthm} for finite groups.   Note that the following
theorem generalizes Theorem \ref{abelian} for finite groups in the
sense that $x$ must be a product of $a$ and $b$.
\begin{theorem}\label{bobsthm}
Suppose $G$ is a finite group and $a,b,x\in G$.  Then
\begin{equation}\label{ineG} Q_p=\Qp \in M_3(\CC)^+ \spa \forall \spa p\in
K(G)\end{equation} if and only if $x=ab$ or $x=ba$.
\end{theorem}
\proof The proof will be outlined in the discussion to follow. The
reverse implication follows immediately from Theorem
\ref{martysthm}, Lemma \ref{kad}, and by noting that we can
identify $K(G)$ with $K(G^{op})$.  The proof in the forward
direction is somewhat lengthy, and the argument is broken into
several cases. The process can be generally described by the
following:

\textbf{Key Steps}
\begin{enumerate}\item Show $x\in\langle a,b \rangle$ (the
subgroup of $G$ generated by $a$ and $b$).  In fact, we
immediately get $x=a^sba^t$ for some natural numbers $s$ and
$t$.\item Show $s+t\equiv 1$ (mod $r$) where $r$ divides $|a|$
(the order of $a$) and is determined by the construction of a
unitary representation of $G$. Thus $x=a^sba^{1-s+qr}$ for some
integer $q$. \item Show $s\equiv 0$ or $s\equiv 1$ (mod $|a|$).
\item Show $q=0$, and thus conclude that $x=ab$ or
$x=ba$.\end{enumerate} We argue the case where $r=1$ and $2\leq
r\leq |a|$ separately.  We treat separately as well the cases
$b\in Hb^{-1}H$ and $b\notin Hb^{-1}H$ where $H=\langle a\rangle$,
the cyclic subgroup of $G$ generated by $a$.

In the next section we include a full proof of Key Step (1) that
illustrates a general approach used in the other steps.  We also
include a proof for the case when $ab=ba^m$ for some $m$ and $b\in
Hb^{-1}H$ to demonstrate the usefulness of the inequalities
derived in the next section.  In virtually all components of the
proofs of the Key Steps above, we repeatedly appeal to these
inequalities.
\section{Fundamental Inequalities}
Assuming that $\{Q_p:p\in K(G)\}\subset M_3(\CC)^+$ leads
immediately to several inequalities involving $a, b$, and $x$,
which are fixed. We define the real valued function $f(p):=$
Det($Q_p$) from $K(G)$ to $\RR$. First, by assumption, $f(p)$ must
be nonnegative, thus
\begin{equation}\label{fund1}
f(p)=1 + 2\real[p(a)p(b)\overline{p(x)}] -
|p(x)|^2-|p(b)|^2-|p(a)|^2 \geq 0 \spa \forall \spa p \in K(G).
\end{equation}  We have already seen the inequality \eqref{fund1}
rearranged to the aesthetically pleasing \eqref{fundin}.  We also
noted that $f(p)=0$ when $p$ is an eigenstate of $a$, and thus by
the assumption that $Q_p\in M_3(\CC)^+$ for all $ p\in K(G)$, such
eigenstates are minima of $f$.

We can always construct positive definite functions on $G$ by
finding a unitary representation $\pi$ of $G$ in $\mnc$ for some
$n$ and applying vector states, i.e., states of the form
$p(\cdot)=\langle \pi (\cdot)\xi,\xi\rangle$ for some unit vector
$\xi\in\CC^n$. In fact, all positive definite functions on $G$
arise in this manner \cite{D}, $\S$13. Also, in this setting,
eigenstates of $a$ are just vector states corresponding to
eigenvectors of $\pi(a)$. For generality, we develop a family of
fundamental inequalities that hold in a $B(\HH)$ setting.

\begin{theorem}\label{nashgen}  Let $a,b,x\in B(\HH)$ where
the operator $a$ is unitary with unit eigenvector $\xi\in\HH$ and
corresponding eigenvalue $\lam$. Suppose $Q_p\in M_3(\CC)^+$ for
all vector states $p$ on $B(\HH)$.  Then for any unit vector
$\eta\in\xi^{\bot}$,
\begin{equation}\label{nashingen}
 |\langle (x-ab)\eta,\xi\rangle + \langle (x-ba)\xi,\eta\rangle |^2
\leq (1-|\langle b\xi,\xi\rangle |^2)\|(\lam I-a)\eta\|^2.
\end{equation}
\end{theorem}
\pf Consider the smooth path $\xi_t=\cos(t)\xi + \sin(t)\eta$
through the unit sphere of $\HH$ parameterized by $t\in\RR$.
Denoting the derivative of $\xi_t$ by $\xi_t'$, note that
$\xi_0=\xi$ and $\xi_0'=\eta$.
 Also, $\xi_t$ induces a smooth path through the vector state space
 of $B(\HH)$ where $p_t$ is
just the vector state determined by $\xi_t$. Thus the function
$f(p_t)=$ Det$(Q_{p_t})$ is a smooth map from $\RR$ to $\RR$, and
classical calculus can be applied. Defining $\hat{a}$ as a
function on $\RR$ by $\hat{a}(t)=p_t(a)$ (similarly for $\hat{x}$
and $\hat{b}$) we can view $f$ as:
\begin{equation}\label{f(p)}
f(t)=1+2\real[\hat{a}(t)\hat{b}(t)\overline{\hat{x}(t)}]-
|\hat{a}(t)|^2 -|\hat{b}(t)|^2-|\hat{x}(t)|^2. \end{equation}

Since $f(0)=0$ is an assumed minimum, we have that $f''(0)\geq 0$.
Recall that $\ah(0)=\lam$ and that $\xh(0)=\lam\bh(0)$ since
$p(x)=p(a)p(b)$ at eigenstates.  Therefore, twice differentiating
$f$ leads to, after some manipulation
\begin{equation}\label{2ndderiv}
|\xh'(0)-\lam\bh'(0)|^2 \leq (|\bh(0)|^2-1)\real
[\overline{\lam}\ah''(0)].
\end{equation}
Applying Leibniz' rule to evaluate \eqref{2ndderiv}, we arrive at
\eqref{nashingen}.\QED

\medskip
We recognize that the family of inequalities captured in
\eqref{nashingen} can be sharpened.
\begin{corollary}\label{sharpnash}
Let $a,b,x \in B(\HH)$ satisfy the hypotheses of Lemma
\ref{nashgen}.  Then for any vector $\eta \in \HH$,
\begin{equation}\label{sharpnashingen}
 |\langle (x-ab)\eta,\xi\rangle| + |\langle (x-ba)\xi,\eta\rangle|
\leq \sqrt{1-|\langle b\xi,\xi\rangle |^2}\|(\lam I-a)\eta\|.
\end{equation}
\end{corollary}
\pf We can choose any $\eta\in\HH$ by writing $\eta=\tilde{\eta} +
c\xi$ for $\tilde{\eta}\in\xi^{\bot}$ and $c\in\CC$ and applying
the unit vector $\frac{e^{i\theta}\eta}{\|\eta\|}$ to
\eqref{nashingen} for an appropriate choice of $\theta$ to arrive
at \eqref{sharpnashingen}. \QED

\medskip

 Of particular interest to us is the case when $B(\HH)= \mnc$, since
 finite groups always have finite dimensional representations as complex
 matrices.
 In particular,
the following corollary holds when the elements $a,b,$ and $x$ of
$G$ are represented as unitary matrices.  In this case we can
always write the matrices with respect to a basis such that $a$ is
diagonal.
\begin{corollary}\label{nash}
Suppose $a,b,x \in \mnc$, where $a$ is diagonal and unitary. We
write $a=$ diag($\lam_1,...,\lam_n$), $b=[b_{jk}]$, and
$x=[x_{jk}]$. Furthermore, suppose $Q_p\in M_3(\CC)^+$ for each
vector state $p$ on $\mnc$.  Then for each $j\in\{1,...,n\}$ we
have
\begin{equation}\label{nashin}
|\sum\limits_{k=1}^n r_k (\xab)|^2\leq(1-|b_{jj}|^2)
\sum\limits_{k=1}^n r_k^2|\lj-\lk|^2 \end{equation} where
\eqref{nashin} holds for any choice of $r_k,\phi_k\in\RR$, for
$k=1,...,n$.
\end{corollary}
\pf We arrive immediately at \eqref{nashin} by applying Lemma
\ref{nashgen} where $\xi = e_j$, the $j^{th}$ standard elementary
basis vector for $\CC^n$, and $\eta$ is defined componentwise by
$\eta_k = r_ke^{i\phi_k}$.  Note that $\eta$ may be arbitrary by
Corollary \ref{sharpnash}.\QED

The family of inequalities captured in \eqref{nashin} are appealed
to repeatedly in the proof of Theorem \ref{bobsthm}.  Although
various other inequalities have been derived pertaining to our
problem, c.f., \cite{C} Ch. 4, these have thus far been the most
useful for our purposes.  In the case of $\mnc$, the corresponding
diagonal entries of $x$, $ab$, and $ba$ when $a$ is diagonal led
us to look for constraints on the off-diagonal elements of $x$.
The above inequalities clearly contain comparisons of $x$ with
$ab$ and $ba$ on the left hand side of \eqref{nashingen}.  The
terms $x_{jk}-\lam_jb_{jk}$ and $x_{kj}-\lam_jb_{kj}$ in
\eqref{nashin} are, when $j\neq k$, the off-diagonal entries from
the matrices $x-ab$ and $x-ba$, respectively (when $a$ is a
diagonal unitary matrix).
\section{Utilizing a Cayley Representation}
Throughout this section we assume $G$ is a finite group of order
$|G|$, and that $a,b,x\in G$.  Also, we use $n_a$ to denote the
order of $a$, i.e., the smallest natural number satisfying
$a^{n_a}=e$,  and $\langle a\rangle$ to denote the cyclic subgroup
of $G$ generated by $a$.  For the majority of the proof of Theorem
\ref{bobsthm} we make use of a Cayley representation $\pi$ of $G$.
\begin{definition}\label{cayleydef}
If $G$ is a finite group enumerated as $g_1,...,g_{|G|}$ we form
the matrix [$g_j^{-1}g_k]\in M_{|G|}(\cs (G))$. This matrix is
sometimes called a \textit{multiplication table} for $G$. Then
associated with this enumeration we have a \textit{Cayley
representation} $\pi$ of $G$ defined by
\[ \pi(g)=[g_{jk}]\mbox{, where }g_{jk}=\left \{ \begin{array}{ll}
1, & \mbox{if }
g=g_j^{-1}g_k \\\\
0, & \mbox{if }  g\neq g_j^{-1}g_k. \end{array}\right.
\]
\end{definition}

Thus, for each $g \in G$, $\pi(g)$ is a $|G| \times |G|$
permutation matrix, hence unitary.  It is easily verified that
$\pi$ is a representation of $G$. Also, $\pi$ is clearly faithful,
and any vector state on $M_{|G|}(\CC)$ gives us a state on $G$.
Another important feature of Cayley representations for our
purposes is highlighted in the next lemma.

\begin{lemma}\label{cayleylem}
Suppose G is a finite group with Cayley representation $\pi$, and
$g,h\in G$. If $\pi(g)$ agrees with $\pi(h)$ in any nonzero entry
then $g=h$.\end{lemma}

\proof This follows directly from the definition of a Cayley
representation. \QED

\medskip
With the generality of a Cayley representation and the sharpness
of Lemma \ref{cayleylem} in mind, we consider a strategic way to
enumerate the group elements.  The general approach will be
described in the proof of the following lemma, which is a proof of
Key Step (1) above.
\begin{lemma}\label{xword}
Suppose $G$ is a finite group and $a,b,x\in G$ satisfy $Q_p \in
M_3(\CC)^+$ for all $p \in K(G)$. Then $x=a^sba^t$ for some
$s,t\in \ZZ_{n_a}$.\end{lemma} \proof If $b\in \langle a\rangle$
or $x\in \langle a \rangle$ then the lemma follows by noting that
any positive definite function $p$ on a subgroup $H$ of $G$ can be
extended to a positive definite function on $G$ by simply letting
$p(g)=0$ for $g\notin H$.  If $p$ is any character on $H=\langle
a\rangle$ extended in this fashion, then from \eqref{fundin} it
follows that $p(x)=p(a)p(b)\neq 0$, implying that both $b\in
\langle a\rangle$ and $x\in \langle a \rangle$.  Therefore
$p(x)=p(ab)$ for each character on $\langle a \rangle$, and hence
$x=ab$. Also, if $a=e$ the lemma follows trivially from
\eqref{fundin}.

 Thus assume $a\neq e$, $b\notin \langle
a\rangle$, and $x\notin \langle a \rangle$. Let $H=\langle
a\rangle$ and choose $\{\beta_1,...,\beta_d\}\subset G$ with
$\beta_1=e$ and $\beta_2=b$, such that $\bigcup_{k=1}^d\beta_kH$
is a coset decomposition of $G$. We enumerate $G$ in such a way
that the top row of the multiplication table reads
\[
e,a,a^2,...,a^{-1},b,ba,ba^2,...,ba^{-1}, \beta_3,\beta_3a,...,
\beta_3a^{-1},......,\beta_d,\beta_da,...,\beta_da^{-1}.
\]
or more concisely, if $H$ is enumerated as $h_1=e, h_2=a,
h_3=a^2$, etc. we may write the top row in shorthand as
\[ H,bH,\beta_3H,...,\beta_dH.\]

We define our Cayley representation $\pi$ by this enumeration of
$G$, and consider $\pi(a)$, $\pi(b)$, and $\pi(x)$ as $|G|$ by
$|G|$ permutation matrices comprised of blocks of size $n_a$ by
$n_a$. Since $\pi$ is faithful we will, for brevity, identify $a$
with $\pi(a)$, etc., without ambiguity. We use the notation
$A_{jk}, B_{jk},$ and $X_{jk}$ to denote the $j,k^{th}$ block of
$a, b$, and $x$ respectively.  In this case, $a$ is completely
determined as a block diagonal matrix where each diagonal block
has the form
\[ A_{jj}=\begin{pmatrix} 0 & 1 & 0 & 0 & ... & 0\\
                               0 & 0 & 1 & 0 & ... & 0\\
                               0 & 0 & 0 & 1 & ... & 0\\
                               . & . & . & . & ... & 0\\
                               . & . & . & . & ... & 1\\
                               1 & 0 & 0 & 0 & ....& 0\\
\end{pmatrix}.\]
With this nice form it is easy to construct eigenvectors (and thus
eigenstates) for $a$.  Also, by assumption, $b$ and $x$ are
matrices with all zeros in the diagonal blocks.  We are led to
consider what form $b$ and $x$ have on off-diagonal blocks, most
importantly those in the first row and column of blocks.  For the
argument in this proof, we need only the upper two block by two
block submatrices of $a,b,$ and $x$, but it demonstrates a general
construction and approach. We use the notation $\#(B_{jk})$ and
$\#(X_{jk})$ to denote the number of nonzero entries in $B_{j,k}$
and $X_{j,k}$, respectively.

 For any $\phi\in\RR$ define the eigenvector $\xi\in \CC^{|G|}$ of $a$
 componentwise by
\[ \xi_j=\left \{ \begin{array}{ll}
\frac{e^{i\phi}}{\sqrt{2n_a}}, & j=1,...,n_a
 \\
\frac{1}{\sqrt{2n_a}}, & j=n_a+1,...,2n_a \\
0, & j=2n_a+1,...,|G|. \end{array}\right.
\]
Let $p$ be the eigenstate of $a$ of the form $p(\cdot)=\langle
\cdot \xi,\xi\rangle$.  Since $p(a)=1$, by \eqref{fundin} we have
$p(x)=p(a)p(b)=p(b)$.  By our construction this translates to
\[\frac{1}{2n_a}[\#(X_{12})e^{-i\phi}
+\#(X_{21})e^{i\phi}]=p(x)=p(b)=\frac{1}{2n_a}[\#(B_{12})e^{-i\phi}
+\#(B_{21})e^{i\phi}]
\]
for any choice of $\phi$.  Therefore, by the linear independence
of $e^{i\phi}$ and $e^{-i\phi}$,  $\#(B_{12})=\#(X_{12})$ and
$\#(B_{21})=\#(X_{21})$.  Since there is at least one nonzero
entry in $B_{12}$, e.g., $b_{1,n_a+1}=1$, there is at least one
nonzero entry in $X_{12}$.  By our construction of the Cayley
representation, this nonzero entry corresponds to an element of
the double coset $HbH$.   Thus, by Lemma \ref{cayleylem},
$x=a^sba^t$ for some $s,t\in \ZZ_{n_a}$. \QED

\medskip
As we go on to Key Steps (2) through (4), the arguments proceed in
similar fashion, beginning with a carefully chosen enumeration of
$G$ and construction of the corresponding Cayley representation.
We choose the enumeration so that $a$ has the same block diagonal
form as in the preceding proof and that $b$ has an easily
manipulatable form.  Next, eigenvectors of $a$ are constructed and
utilized along with Lemma \ref{nashin} to further restrict the
form of $x$ as a word in $a$ and $b$.  We take full advantage of
the fact that the eigenvalues of $a$ are the $n_a^{th}$ roots of
unity, which work nicely in the context of Lemma \ref{nashin}.

As the steps for the different cases proceed the choice of
enumeration and construction of eigenvectors becomes unexpectedly
more delicate and involved, see \cite{C}, but follow in the same
spirit as described in the proof of Lemma \ref{xword}.  We include
a proof for the particular case when the relationship $ab=ba^m$
holds for some $m$ and $b\in H b^{-1} H$. An example of a finite
group for which this case pertains would be any of the quaternion
groups.  The purpose of including this proof is to give the reader
a small sample of the usefulness of Theorem \ref{nashgen} and its
corollaries.
\begin{theorem}\label{semidirect}
Suppose $G$ is a finite group and $a,b,x\in G$ satisfy $Q_p \in
M_3(\CC)^+$ for all $p\in K(G)$.  Suppose further that $ab=ba^m$
for some $m\in\ZZ_{n_a}$ and that $b\in Hb^{-1}H$. Then $x=ab$ or
$x=ba$.
\end{theorem}
\pf If $a=e$, $b\in \langle a\rangle$, or $x\in \langle a \rangle$
then $x=ab=ba$ by the first paragraph of the proof of Lemma
\ref{xword}. Thus assume $a\neq e$, $b\notin \langle a\rangle$,
and $x\notin \langle a \rangle$.

We know by Lemma \ref{xword} that $x=a^sba^t$ for some $s$ and
$t$.  Therefore, since $ab=ba^m$, we have that $x=ba^l$ for some
$l\in\ZZ_{n_a}$.  We begin by showing  that $l\equiv 1+q(m-1)$
(mod $n_a$) for some $q\in\ZZ_{n_a}$.

 Let $\pi$ be the same Cayley representation of $G$
used in the proof of Lemma \ref{xword} and suppose
 $\lam=e^{2\pi
i/(m-1,n_a)}$, where $(m-1,n_a)$ denotes the greatest common
divisor of $m-1$ and $n_a$. Again, for brevity, identify $a$ with
$\pi(a)$, etc. Let the eigenvector $\xi$ of $a$ be defined
componentwise by
\begin{equation}\label{xi}
\xi_j=\left \{
\begin{array}{ll} \frac{\lam^j}{\sqrt{2n_a}}e^{i\phi}, &
j=1,2,...,n_a
 \\
 \frac{\lam^j}{\sqrt{2n_a}}, & j=n_a+1,n_a+2...,2n_a
 \\
0, & \mbox{otherwise.} \end{array}\right.
\end{equation}
where $\phi$ is arbitrary.   Then $\xi$ is a unit eigenvector of
$a$ with eigenvalue $\lam$.   It follows, since $ab=ba^m$, that
$b\xi$ is an eigenvector of $a$ with eigenvalue $\lam^m=\lam$. By
our constructions it follows that, for some $z\in\CC$, we have
$\langle b\xi,
\xi\rangle=\frac{1}{2n_a}(n_ae^{-i\phi}+ze^{i\phi})$.  We choose
$\phi$ so that $\langle b\xi,\xi\rangle \neq 0$.

Recall that $x$ must agree with $ab$ at eigenstates of $a$. Thus
$\lam\langle b\xi,\xi\rangle = \langle ab\xi,\xi\rangle =\langle
x\xi, \xi\rangle = \langle ba^l\xi,\xi\rangle =\lam^l\langle
b\xi,\xi\rangle$. In other words, we have $\lam^l=\lam$ and hence
$l-1\equiv 0$ (mod $(m-1,n_a)$).  Therefore, $l=1+s(m-1,n_a)$ for
some $s\in\ZZ$, giving $l\equiv 1+q(m-1)$ (mod $n_a$) for some
$q\in \ZZ_{n_a}$, since $(m-1,n_a)$ is congruent, mod $n_a$, to a
multiple of $m-1$.

 For uniqueness, we
may assume $q\in \{0,1,...,|m-1|-1\}$ (where $|m-1|$ denotes the
order of $(m-1)$ in $\ZZ_{n_a}$).  The next step is to show that
$q$ must be 0 or 1.

Note that $ab=ba^m$ implies $a^kb=ba^{km}$ for each $k\in\NN$ and
hence $k\equiv 0$ (mod $n_a$) if and only if $km \equiv 0$ (mod
$n_a$). Thus $m$ has order $n_a$ in $\ZZ_{n_a}$, and hence there
is some unique $m^{-1}\in\ZZ_{n_a}$ such that $mm^{-1}\equiv 1$
(mod $n_a$).  It follows that $x=ba^l=a^{m^{-1}l}b$.  Let $\lam$
be an arbitrary eigenvalue of $a$ and define the eigenvector $\xi$
of $a$ as in \eqref{xi}. Letting $\eta =b\xi$ and writing
$x=a^{m^{-1}l}b$ in the first term and $x=ba^l$ in the second term
of \eqref{sharpnashingen}, we have
\[ |\lam^{m^{-1}l}-\lam||\langle b^2\xi,\xi\rangle|+|\lam^l-\lam|
\leq |\lam-\lam^m|.\] However, since $b\in Hb^{-1}H$ it follows
that $\langle b^2\xi,\xi\rangle=\lam^t$ for some $t\in \ZZ_{n_a}$.
Also, since $l=1+q(m-1)$, the previous inequality may be written
\[ |\lam^{(1-q)(m^{-1}-1)}-1|+|\lam^{q(m-1)}-1|
\leq |\lam^{m-1}-1|.\]

Using the fact that $\{\lam^{m-1}|\lam$ is an eigenvalue of $a\}$
is a subgroup of the $n_a^{th}$ roots of unity, we select $\lam$
so that $\lam^{m-1}$ is as close as possible (but not equal) to
$1$. Specifically we choose $\lam$ such that $\lam^{m-1}=e^{2\pi
i/|m-1|}$.  Noting that $|\lam^{q(m-1)}-1| \leq |\lam^{m-1}-1|$ we
can see that the only choices for $q$ are $0$, $1$, and $|m-1|-1$.
If $q=|m-1|-1$ then we must have $\lam^{2(m^{-1}-1)}=1$ for all
eigenvalues $\lam$ of $a$.  Since $|m-1|=|m^{-1}-1|$, this implies
that $2=|m^{-1}-1|=|m-1|$, and thus the only choices for $q$ are
$0$ or $1$.  Therefore, $q=0$ or $q=1$, implying $x=ba$ or
$x=ba^m=ab$.

\QED

\medskip
 Again, the previous proof only gives a glimpse of the technique
used for each of the various cases. The reader may note that the
proof needs only minor modifications to obtain the case where
$b\notin Hb^{-1}H$.  Refined versions of the proofs in \cite{C}
will appear shortly.

\bibliographystyle{amsalpha}

\end{document}